\documentclass[11pt]{amsart}

\usepackage{amssymb,epstopdf,subcaption,mathtools,tikz-cd,pb-diagram,thmtools, thm-restate,pinlabel,pifont,enumitem,stmaryrd}
\usepackage{graphicx}

\usepackage{tikz}
\usepackage{dsfont}
\usetikzlibrary{matrix}

\usepackage[all]{xy}
\usepackage[mathscr]{euscript}

\newcommand{\nc}{\newcommand}
\nc{\dmo}{\DeclareMathOperator}
\nc{\nt}{\newtheorem}

\newtheorem{theorem}{Theorem}[section]
\newtheorem{mainthm}{Theorem}

\newtheorem{proposition}[theorem]{Proposition}
\newtheorem{corollary}[theorem]{Corollary}
\theoremstyle{definition}

\newtheorem{question}[theorem]{Question}

\newtheorem{remark}[theorem]{Remark}
\theoremstyle{remark}

\nc{\cut}{\!\ssearrow\!}

\dmo{\Diff}{Diff}
\dmo{\Mod}{Mod}
\dmo{\SMod}{SMod}
\dmo{\I}{\mathcal{I}}
\dmo{\SO}{SO}
\dmo{\Orth}{O}
\dmo{\Sp}{Sp}
\dmo{\SL}{SL}
\dmo{\GL}{GL}
\dmo{\im}{im}
\dmo{\Emb}{Emb}
\dmo{\PSp}{PSp}
\dmo{\PSL}{PSL}
\dmo{\PMod}{PMod}
\dmo{\Homeo}{Homeo}
\dmo{\Twist}{Twist}
\dmo{\Aut}{Aut}
\dmo{\Nil}{Nil}
\dmo{\Sol}{Sol}
\dmo{\Isom}{Isom}
\dmo{\Out}{Out}
\dmo{\OUT}{\mathcal{OUT}}
\dmo{\AUT}{\mathcal{AUT}}
\dmo{\Inn}{Inn}
\dmo{\ab}{ab}
\dmo{\orb}{orb}
\dmo{\Hom}{Hom}
\dmo{\Push}{Push}

\nc{\Z}{\mathbb Z}
\nc{\N}{\mathcal N}
\nc{\R}{\mathbb R}
\nc{\F}{\mathcal F}
\nc{\C}{\mathbb{C}}

\nc{\ga}{\gamma}
\nc{\de}{\delta}
\nc{\ep}{\epsilon}

\nc{\flm}{\lambda_{2}}

\nc{\normalclosure}[1]{\ensuremath{\left \langle \left \langle #1 \right \rangle \right \rangle}}

\nc{\margin}[1]{\marginpar{\scriptsize #1}}

\nc{\p}[1]{\bigskip\noindent\textbf{#1.}}
\nc{\lei}[1]{{\color{red} \sf  L: [#1]}}
\nc{\bena}[1]{{\color{blue} \sf  B: [#1]}}

\title{Mapping class groups of circle bundles over a surface}
\author{Lei Chen}
\author{Bena Tshishiku}

\address{Lei Chen \\ Department of Mathematics\\ University of Maryland \\ 4176 Campus Drive\\ College Park, MD 20742, USA \\  chenlei@umd.edu}

\address{Bena Tshishiku \\ Department of Mathematics\\ Brown University \\ 151 Thayer St.  \\ Providence, RI, 02912, USA\\  bena$\_$tshishiku@brown.edu.   }

\date{\today}

\begin{document}

\maketitle

\vspace*{-4ex}

\begin{abstract}
In this paper, we study the algebraic structure of mapping class group $\Mod(X)$ of 3-manifolds $X$ that fiber as a circle bundle over a surface $S^1\to X\to S_g$. There is an exact sequence $1\to H^1(S_g)\to \Mod(X)\to\Mod(S_g)\to1$. We relate this to the Birman exact sequence and determine when this sequence splits. 
\end{abstract}

\section{Introduction}

For $g\ge1$, let $S_g$ denote the closed oriented surface of genus $g$, and for $k\in\Z$, let $X_{g}^{k}$ denote the closed 3-manifold that fibers 
\[S^1\to X_{g}^{k}\to S_g\]
as an oriented circle-bundle with Euler number $k$. Assuming $(g,k)\neq(1,0)$, the mapping class group $\Mod(X_{g}^{k}):=\pi_0\big(\Homeo^+(X_{g}^{k})\big)$ fits into a short exact sequence 
\begin{equation}\label{eqn:Mod-SES}1\rightarrow H^1(S_g;\Z)\rightarrow\Mod(X_{g}^{k})\rightarrow\Mod(S_g)\rightarrow1.\end{equation}

This paper is motivated by the following question. 

\begin{question}\label{q:main}
For which values of $g,k$ is the extension in (\ref{eqn:Mod-SES}) split? 
\end{question}

Interestingly, the extension does split for $k=2-2g$, in which case $X_g^{k}$ is unit tangent bundle $US_g$. In fact, there is a natural action of $\Mod(S_g)$ on $US_g$ by homeomorphisms, which gives a splitting of (\ref{eqn:Mod-SES}) upon taking isotopy classes. For $g\ge2$, this action comes from the action of the punctured mapping class group $\Mod(S_{g,1})$ on triples of points on the boundary of hyperbolic space $\mathbb H^2$. This construction dates back to the work of Nielsen. See \cite[\S5.5.4, \S8.2.6]{Farb-Margalit} and \cite[\S1]{Souto}. 

In general, Question \ref{q:main} reduces to a question about group cohomology. 
The extension (\ref{eqn:Mod-SES}) splits if and only if its Euler class $eu_k\in H^2\big(\Mod(S_g);H^1(S_g;\Z)\big)$ vanishes \cite[\S IV.3]{Brown}. Here the coefficients are twisted via the natural action of $\Mod(S_g)$ on $H^1(S_g;\Z)$. However, a computation of $H^2\big(\Mod(S_g);H^1(S_g;\Z)\big)$ does not appear to be in the literature. 

The extension \eqref{eqn:Mod-SES} is related to the Birman exact sequence
\[1\to\pi_1(S_g)\to\Mod(S_{g,1})\to\Mod(S_g)\to1.\]
By taking quotients by the commutator subgroup $\pi'\equiv [\pi_1(S_g),\pi_1(S_g)]$, we obtain the following extension
\begin{equation}\label{eqn:BES-quotient}1\to H_1(S_g)\to\Mod(S_{g,1})/\pi'\to\Mod(S_g)\to1.\end{equation}
Our main result relates the sequences \eqref{eqn:Mod-SES} and \eqref{eqn:BES-quotient}.

\begin{mainthm}\label{thm:main}
Fix $g\ge1$ and $k\in\Z$. Assume $(g,k)\neq(1,0)$. There is a map between the short exact sequences \eqref{eqn:Mod-SES} and \eqref{eqn:BES-quotient}
\[\begin{xy}
(-60,0)*+{1}="A";
(-30,0)*+{H_1(S_g)}="B";
(0,0)*+{\Mod(S_{g,1})/\pi'}="C";
(30,0)*+{\Mod(S_g)}="D";
(60,0)*+{1}="E";
(-60,-15)*+{1}="F";
(-30,-15)*+{H^1(S_g;\Z)}="G";
(0,-15)*+{\Mod(X_{g}^{k})}="H";
(30,-15)*+{\Mod(S_g)}="I";
(60,-15)*+{1}="J";
{\ar "A";"B"}?*!/_3mm/{};
{\ar "B";"C"}?*!/_3mm/{};
{\ar "C";"D"}?*!/^3mm/{};
{\ar "D";"E"}?*!/^3mm/{};
{\ar "F";"G"}?*!/_3mm/{};
{\ar "G";"H"}?*!/_3mm/{};
{\ar "H";"I"}?*!/^3mm/{};
{\ar "I";"J"}?*!/^3mm/{};
{\ar "B";"G"}?*!/_3mm/{k\delta};
{\ar "C";"H"}?*!/^3mm/{};
{\ar@{=} "D";"I"}?*!/_3mm/{}
\end{xy}\]
The homomorphism $k\delta$ is the Poincar\'e duality isomorphism $\delta$ composed with multiplication by $k$. In particular, when $k=1$, the exact sequences \eqref{eqn:Mod-SES} and \eqref{eqn:BES-quotient} are isomorphic. 
\end{mainthm}

Theorem \ref{thm:main} implies the Euler classes of the extensions (\ref{eqn:Mod-SES}) satisfy $eu_k=k\>eu_1$ for fixed $g$. Next we determine the subgroup generated by $eu_1$ in $H^2\big(\Mod(S_g);H^1(S_g;\Z)\big)$. 

\begin{theorem}\label{thm:H2H1}
Fix $g\ge1$, and let $eu_1$ be the Euler class of the extension \eqref{eqn:Mod-SES}. Then $eu_1$ has order $2g-2$ in $H^2\big(\Mod(S_g);H_1(S_g;\Z)\big)$. Furthermore, if $g\ge 8$, then $eu_1$ generates this group, i.e.\
\[H^2\big(\Mod(S_g);H^1(S_g;\Z)\big)\cong \mathbb{Z}/(2g-2)\Z.\] 
\end{theorem}

Combining Theorem \ref{thm:main} and Theorem \ref{thm:H2H1} we obtain the following answer to Question \ref{q:main}. 


\begin{corollary}\label{cor:splitting}
For $g\ge2$ and $k\in\mathbb Z$, the extension \eqref{eqn:Mod-SES} splits if and only if $k$ is divisible $2g-2$. For $g=1$ the extension splits for each $k$. 
\end{corollary}

When a splitting exists, the different possible splittings (up to the action of $H^1(S_g;\Z)$ on $\Mod(X_{g}^{k})$ by conjugation) are parameterized by elements of  $H^1\big(\Mod(S_g); H^1(S_g;\Z)\big)$ \cite[Ch.\ IV, Prop.\ 2.3]{Brown}. This group vanishes for $g\ge1$ \cite[Prop.\ 4.1]{Morita}, so the splitting, when it exists, is unique.

\p{Connection to Nielsen realization}
Instead of Question \ref{q:main}, one can ask whether there is a splitting of the composite surjection
\[\Homeo(X_{g}^{k}) \to\Mod(X_{g}^{k})\to\Mod(S_g).\] This is an instance of a Nielsen realization problem. Of course, if $\Mod(X_{g}^{k})\to\Mod(S_g)$ does not split, then neither does $\Homeo(X_{g}^{k})\to\Mod(S_g)$, and Corollary \ref{cor:splitting} gives examples of this. Since $\Mod(S_g)$ has a natural action on $US_g$, the surjection $\Homeo(X_{g}^{k})\to\Mod(S_g)$ does split for $k=\pm(2g-2)$. This is somewhat surprising since mapping class groups are rarely realized as groups of surface homeomorphisms \cite{Markovic, Chen, Chen-Salter}. We wonder whether this splitting is unique, or if a splitting exists for other values $k$ divisible by $2g-2$ (for example, $k=0$). We plan to study this in a future paper.

\p{Previous work and proof techniques} Waldhausen \cite[\S7]{Waldhausen} proved that the group $\pi_0\big(\Homeo(X_{g}^{k})\big)$ is isomorphic to the outer automorphism group $\Out\big(\pi_1(X_{g}^{k})\big)$. 
From this, the short exact sequence \eqref{eqn:Mod-SES} can be derived from work of Conner--Raymond \cite{Conner-Raymond} and the Dehn--Nielsen--Baer theorem; alternatively, see McCullough \cite[\S3]{McCullough-virtual}. The Dehn--Nielsen--Baer theorem also plays a central role in Theorem \ref{thm:main}, since it allows us to translate back and forth between topology and group theory. There is a mix of both in the proof of Theorem \ref{thm:main} in \S\ref{sec:main}. 

To prove Theorem \ref{thm:main}, we consider a version of Question \ref{q:main} where $X_g^k$ and $S_g$ are punctured. For the punctured manifolds, similar to \eqref{eqn:Mod-SES}, there is a short exact sequence 
\[1\to H^1(S_g;\Z)\to\Mod(X_{g,1}^k)\to\Mod(S_{g,1})\to 1,\]
and we construct a splitting 
\[\sigma:\Mod(S_{g,1})\to\Mod(X_{g,1}^k).\] See Corollary \ref{cor:Aut-SES}. 
A key part of our proof of Theorem \ref{thm:main} is to determine the image of the point-pushing subgroup $\pi_1(S_g)<\Mod(S_{g,1})$ under $\sigma$. For this we relate three natural surface group representations $\pi_1(S_g)\to\Mod(X_{g,1}^k)$ that appear in the following diagram, where the diagonal map is point pushing on $X_g^k$ (not a commutative diagram). 
\[\begin{xy}
(-30,0)*+{\pi_1(S_g)}="A";
(40,0)*+{\Mod(S_{g,1})}="B";
(-30,-15)*+{H^1(S_g;\Z)}="C";
(40,-15)*+{\Mod(X_{g,1}^k)}="D";
(5,-8)*+{}="E";
{\ar "A";"B"}?*!/_3mm/{\text{point-pushing on $S_g$}};
{\ar "B";"D"}?*!/_3mm/{\sigma};
{\ar "A";"C"}?*!/^3mm/{};
{\ar "A";"D"}?*!/^3mm/{};
{\ar "C";"D"}?*!/^2mm/{\text{transvections}}
\end{xy}\]
See Proposition \ref{prop:section-on-inner} for a precise statement.


In order to deduce Corollary \ref{cor:splitting}, we use a spectral sequence argument to prove that $eu_1$ generates a subgroup of $H^2\big(\Mod(S_g);H^1(S_g;\Z)\big)$ isomorphic to $\Z/(2g-2)\Z$. A different spectral sequence computation proves that $eu_1$ generates $H^2\big(\Mod(S_g);H^1(S_g;\Z)\big)$ when $g$ is large. These computations use several known computations, including work of Morita \cite{Morita}.

\p{Section outline} In \S\ref{sec:background} we collect the results we need about the manifolds $X_{g}^{k}$ and their mapping class groups, including Waldhausen's work. Theorem \ref{thm:main} is proved in \S\ref{sec:main}; this section is the core of the paper. In \S\ref{sec:sseq}, we do two spectral sequence computations to prove Theorem \ref{thm:H2H1}. 

\p{Acknowledgement} Thanks to B.\ Farb for sharing the reference \cite{McCullough-virtual} and to D.\ Margalit for comments on a draft. The authors are supported by NSF grants DMS-2203178, DMS-2104346 and DMS-2005409. 

\section{Circle bundles over surfaces} \label{sec:background}

Here we review some results about circle bundles over surfaces that we will need in future sections. 

\subsection{Classification}
By an oriented circle bundle we mean a fiber bundle 
\[S^1\to E\to B\]with structure group $\Homeo^+(S^1)$, the group of orientation preserving homeomorphisms of $S^1$. The inclusion of the rotation group $\SO(2)$ in $\Homeo^+(S^1)$ is a homotopy equivalence, so circle bundles are in bijection with rank-2 real vector bundles. The classifying space $B\SO(2)$ is homotopy equivalent to $\C P^\infty$, which is an Eilenberg--Maclane space $K(\Z,2)$. Thus each circle bundle is uniquely determined up to isomorphism by its Euler class $eu(E)\in H^2(B;\Z)$, which is the primary obstruction to a section of the bundle. 

When $B=S_g$ is a closed, oriented surface, $H^2(S_g;\Z)\cong\Z$, we can speak of the Euler \emph{number}.  We use $X_{g}^{k}$ to denote the total space of the circle bundle 
\[S^1\to X_{g}^{k}\to S_g\] with Euler number $k$. For example, for the unit tangent bundle $eu(US_g)=2-2g$ (the Euler characteristic), so $US_g\cong X_{g}^{2-2g}$. We also note that $X_{g}^{k}$ and $X_{g}^{-k}$ are homeomorphic 3-manifolds, since the sign of the Euler number of a circle bundle over $S_g$ depends on the choice of orientation on $S_g$. 

\subsection{\boldmath Fundamental group $\pi_1(X_{g}^{k})$ and its automorphisms}\label{sec:pi1} From the long exact sequence of a fibration, we have an exact sequence
\[1\to\Z\to\pi_1(X_{g}^{k})\to\pi_1(S_g)\to1.\]

The group $\pi_1(X_{g}^{k})$ has a presentation 
\begin{equation}\label{eqn:presentation}\pi_1(X_{g}^{k})=\big\langle A_1,B_1,\ldots,A_g,B_g,z\mid z\text{ central, } [A_1,B_1]\cdots[A_g,B_g]=z^k
\big\rangle.\end{equation}
Using this, one finds $\langle z\rangle\cong\Z$ is the center of $\pi_1(X_{g}^{k})$ as long as 
$(g,k)\neq(1,0)$. When $g\ge2$ follows from the fact that the group $\pi_1(S_g)$ has trivial center; the case $g=1$ can be treated directly.

Given this computation of the center, any automorphism of $\pi_1(X_{g}^{k})$ induces an automorphism of $\langle z\rangle\cong\Z$ and descends to an automorphism of $\pi_1(S_g)$. The latter gives a homomorphism
\[\Aut\big(\pi_1(X_{g}^{k})\big)\to\Aut\big(\pi_1(S_g)\big)\]
that restricts to an isomorphism between the inner automorphism groups 
\begin{equation}\label{eqn:inner}\Inn\big(\pi_1(X_{g}^{k})\big)\cong\pi_1(S_g)\cong\Inn\big(\pi_1(S_g)\big)\end{equation}
and hence descends to a homomorphism 
\begin{equation}\label{eqn:Out-map}\Out\big(\pi_1(X_{g}^{k})\big)\to\Out\big(\pi_1(S_g)\big).\end{equation}

\p{Orientations} It will be convenient to define 
\[\AUT\big(\pi_1(S_g)\big)<\Aut\big(\pi_1(S_g)\big)\]
as the subgroup that acts trivially on $H_2(\pi_1(S_g);\Z)\cong\Z$ (the ``orientation-preserving" subgroup). We define 
\[\AUT\big(\pi_1(X_{g}^{k})\big)<\Aut\big(\pi_1(X_{g}^{k})\big)\] as the group of automorphisms that project into $\AUT\big(\pi_1(S_g)\big)$ and that act trivially on the center $\langle z\rangle\cong\Z$. In particular, $\AUT\big(\pi_1(X_{g}^{k})\big)$ has index 4 in $\Aut\big(\pi_1(X_{g}^{k})\big)$. 

These orientation-preserving subgroups contain the (respective) inner automorphism groups, and we denote the quotients $\OUT\big(\pi_1(X_{g}^{k})\big)$ and $\OUT\big(\pi_1(S_g)\big)$. 

\subsection{\boldmath Mapping class group $\Mod(X_{g}^{k})$}\label{sec:Waldhausen}

Fix $g\ge1$ and $k\in\Z$, and assume $(g,k)\neq(1,0)$. Let $\Homeo^+(X_{g}^{k})$ denote the group of homeomorphisms whose image in $\Out\big(\pi_1(X_{g}^{k})\big)$ is contained in $\OUT\big(\pi_1(S_g)\big)$. Define
\[\Mod(X_{g}^{k}):=\pi_0\big(\Homeo^+(X_{g}^{k})\big).\]

Waldhausen \cite[Cor.\ 7.5]{Waldhausen} proved that the natural homomorphism
\[\pi_0\big(\Homeo(X_{g}^{k})\big)\to\Out\big(\pi_1(X_{g}^{k})\big)\]
is an isomorphism. Then, by the definitions, this homomorphism restricts to an isomorphism $\Mod(X_{g}^{k})\cong\OUT\big(\pi_1(X_{g}^{k})\big)$. Waldhausen also proved that 
$\pi_0\Homeo(X_{g}^{k})$ is isomorphic to the group of fiber-preserving homeomorphisms modulo homeomorphisms that are isotopic to the identity through fiber-preserving isotopies; see \cite[Rmk.\ following Cor.\ 7.5]{Waldhausen}. Consequently, there is a homomorphism 
\begin{equation}\label{eqn:Mod-map}\Mod(X_{g}^{k})\to\Mod(S_g).\end{equation}

Altogether, we have the following commutative diagram relating the maps (\ref{eqn:Out-map}) and (\ref{eqn:Mod-map}). 
\begin{equation}\label{eqn:Mod-Out}
\begin{xy}
(-20,0)*+{\Mod(X_{g}^{k})}="A";
(20,0)*+{\Mod(S_g)}="B";
(-20,-15)*+{\OUT\big(\pi_1(X_{g}^{k})\big)}="C";
(20,-15)*+{\OUT\big(\pi_1(S_g)\big)}="D";
{\ar "A";"B"}?*!/_3mm/{};
{\ar "B";"D"}?*!/_3mm/{\cong};
{\ar "A";"C"}?*!/^3mm/{\cong};
{\ar "C";"D"}?*!/_3mm/{}
\end{xy}\end{equation}
The right vertical map is an isomorphism by the Dehn--Nielsen--Baer theorem \cite[Thm.\ 8.1]{Farb-Margalit}. 
Furthermore, by Conner--Raymond \cite[Thm.\ 8]{Conner-Raymond} that there is a short exact sequence 
\begin{equation}\label{eqn:Out-SES}1\to\Hom(\pi_1(S_g),\Z)\to\OUT\big(\pi_1(X_{g}^{k})\big)\to\OUT\big(\pi_1(S_g)\big)\to1.\end{equation}
This establishes the short exact sequence \eqref{eqn:Mod-SES} in the introduction. We will give a concrete derivation of this exact sequence in Corollary \ref{cor:Aut-SES} below.

\section{Relating $\Mod(X_{g}^{k})$ to the Birman exact sequence}\label{sec:main}

In this section, we prove Theorem \ref{thm:main}. To construct the map of short exact sequences in Theorem \ref{thm:main}, our main task is to first define a homomorphism $\Mod(S_{g,1})\to\Mod(X_{g}^{k})$ and then to compute that its kernel is the commutator subgroup of $\pi_1(S_g)<\Mod(S_{g,1})$ (the point-pushing subgroup). We do this is \S\ref{sec:hom} and \S\ref{sec:kernel}.

\subsection{\boldmath A homomorphism $\Psi:\Mod(S_{g,1})\to\Mod(X_{g}^{k})$}\label{sec:hom}

Fix a basepoint $*\in S_g$, and set $S_{g,1}=S_g\setminus\{*\}$. 
By the Dehn--Nielsen--Baer theorem, $\Mod(S_{g,1})$ is isomorphic to $\Out^*(F_{2g})$, where $F_{2g}$ is the free group of rank $2g$ and $\Out^*(F_{2g})<\Out(F_{2g})$ is the subgroup that preserves the conjugacy class corresponding to the free homotopy class of the curve around the puncture in $S_{g,1}$. We construct $\Psi$ as a composition
\begin{equation}\label{eqn:Psi}\Psi:\Mod(S_{g,1})\cong\Out^*(F_{2g})\xrightarrow{\sigma}\AUT\big(\pi_1(X_{g}^{k})\big)\to\OUT\big(\pi_1(X_{g}^{k})\big)\cong\Mod(X_{g}^{k}).\end{equation}

To define $\sigma$, fix a generating set $\alpha_1,\beta_1,\ldots,\alpha_g,\beta_g$ for $F_{2g}$ such that $c=\prod_{i=1}^g[\alpha_i,\beta_i]$ represents the conjugacy class of the curve around the puncture. Let 
\begin{equation}\label{eqn:iota}\iota:F_{2g}\to\pi_1(X_{g}^{k})\end{equation}
 be the homomorphism defined by $\alpha_i\mapsto A_i$ and $\beta_i\mapsto B_i$. Given $f\in\Out^*(F_{2g})$, fix an automorphism $\tilde f:F_{2g}\to F_{2g}$ that represents $f$, and assume that $\tilde f(c)=c$ (this can always be achieved by composing any lift with an inner automorphism of $F_{2g}$). Next we define $\sigma(f)$ on generators of $\pi_1(X_{g}^{k})$ by
\begin{equation}\label{eqn:sigma}\sigma(f)(A_i)=\iota\tilde f(\alpha_i),\>\>\>\sigma(f)(B_i)=\iota\tilde f(\beta_i),\>\>\>\sigma(f)(z)=z.\end{equation}
To show that $\sigma(f)$ extends to a homomorphism of $\pi_1(X_{g}^{k})$, we check that the relation $[A_1,B_1]\cdots[A_g,B_g]=z^k$ is preserved under $\sigma(f)$: 
\[\prod_i[\sigma(f)(A_i),\sigma(f)(B_i)]=\prod_i[\iota\tilde f(\alpha_i),\iota\tilde f(\beta_i)]=\iota(c)=z^k=\sigma(f)(z^k).\]
The second equality uses the the fact that $\tilde f(c)=c$. The map $\sigma(f)$ is independent of the choice of $\tilde f$ because different choices of $\tilde f$ differ by conjugation by powers of $c$ (because the centralizer of $c$ in $F_{2g}$ is the cyclic subgroup $\langle c\rangle$)\footnote{Note that the centralizer is isomorphic to $\Z$ and contains $\langle c\rangle$. It is only bigger if $c=u^i$ for some $u\in F_{2g}$ and $i\ge2$. By contradiction, if $c=u^i$ for $i\ge2$, then $u$ is cyclically reduced because $c$ is. This implies that $u$ is a subword of $c=\prod[\alpha_i,\beta_i]$, which is absurd. } and $\iota(c)=z^k$ is central in $\pi_1(X_{g}^{k})$. 
The homomorphism $\sigma(f):\pi_1(X_{g}^{k})\to\pi_1(X_{g}^{k})$ is an automorphism and belongs to $\AUT\big(\pi_1(X_{g}^{k})\big)$ by definition. Furthermore, $f\mapsto\sigma(f)$ is a homomorphism, which is easy to check using the observation that if $w=\iota w'$, then $\sigma(f)(w)=\iota\tilde f(w')$.

Composing $\sigma$ with $\AUT\to\OUT$ gives the desired homomorphism $\Psi$.  As a corollary of this construction, we have proved the following. 

\begin{corollary}\label{cor:Aut-SES}
Fix $g\ge1$ and $k\in\Z$, and assume $(g,k)\neq(1,0)$. The natural map $\Phi: \AUT\big(\pi_1(X_{g}^{k})\big)\to\AUT\big(\pi_1(S_g)\big)$ (see \S$\ref{sec:pi1}$) fits into an exact sequence
\begin{equation}\label{eqn:Aut-SES}1\to\Hom(\pi_1(S_g),\Z)\to\AUT\big(\pi_1(X_{g}^{k})\big)\xrightarrow{\Phi}\AUT\big(\pi_1(S_g)\big)\to1,
\end{equation}
and this exact sequence splits. 
\end{corollary}

\begin{proof}
First we compute the kernel of $\Phi$. Using the presentation for $\pi_1(X_{g}^{k})$ in (\ref{eqn:presentation}), if $f\in\ker(\Phi)$, then 
\[f(A_i)=A_iz^{m_i}\>\>\>\text{ and }\>\>\>f(B_i)=B_iz^{n_i}\]
for some $m_1,n_1,\ldots,m_g,n_g\in\Z$. The map $a_i\mapsto m_i$, $b_i\mapsto n_i$ extends to a homomorphism $\tau(f):\pi_1(S_g)\to\Z$. It is elementary to check that the map $\ker(\Phi)\to\Hom(\pi_1(S_g),\Z)$ defined by $f\mapsto \tau(f)$ is an isomorphism. 

The homomorphism $\sigma$ defined above shows that $\Phi$ is a split surjection. Note that the $\Mod(S_g,*)\cong\Mod(S_g\setminus\{*\})$ (basepoint vs.\ puncture), so by Dehn--Nielsen--Baer there is an isomorphism $\AUT\big(\pi_1(S_g)\big)\cong\Out^*(F_{2g})$, and we use this isomorphism to view $\sigma$ as a splitting of $\Phi$. 
\end{proof}

\begin{remark}
We call elements of $\ker(\Phi)\cong H^1(S_g;\Z)$ \emph{transvections}. 
\end{remark}

\begin{remark}
The homomorphism $\Psi$ can be constructed on the level of topology as follows. Fix a regular neighborhood $D$ of the puncture on $S_{g,1}$ (so $D$ is a once-punctured disk). Given a mapping class $f\in\Mod(S_{g,1})$, choose a representing homeomorphism $\mathfrak f$. Without loss of generality, we can assume that $\mathfrak f$ is the identity on $D$. The bundle $X_g^k\to S_g$ can be trivialized over $S\setminus D$ (because the classifying space $B\SO(2)$ is simply connected). Fixing a trivialization $(S\setminus D)\times S^1$ over $S\setminus D$, we lift $\mathfrak f$ to the product homeomorphism $\mathfrak f\times\text{id}_{S^1}$. This homeomorphisms is the identity on the boundary $\partial(S\setminus D)\times S^1$, so we can extend by the identity to obtain a homeomorphism $\tilde{\mathfrak f}$ of $X_g^k$. The map sending $f\in\Mod(S_{g,1})$ to the isotopy class $\left[\tilde{\mathfrak f}\right]\in\Mod(X_g^k)$ is the topological version of the homomorphism $\Psi$. Note that the isotopy class $[\mathfrak f]$ is only well-defined up to Dehn twists about $\partial D$ which is a loop around the puncture. This is analogous to the ambiguity encountered in the definition of $\sigma$, which ultimately does not affect the definition of $\Psi$. 
\end{remark}

Corollary \ref{cor:Aut-SES} and equation \eqref{eqn:inner} combine to give the short exact sequence of outer automorphism groups \eqref{eqn:Out-SES}. 

\p{Warning} The splitting of the short exact sequence (\ref{eqn:Aut-SES}) does \emph{not} give a splitting of the short exact sequence (\ref{eqn:Out-SES}). Indeed we will show the latter sequence does \emph{not} always split (Corollary \ref{cor:splitting}). The subtlety comes from the fact that the inner automorphism group $\Inn\big(\pi_1(X_{g}^{k})\big)\cong\pi_1(S_g)$ does not coincide with the image of $\pi_1(S_g)\cong\Inn\big(\pi_1(S_g)\big)<\Aut\big(\pi_1(S_g)\big)$ under the section $\sigma$. Proposition \ref{prop:section-on-inner} below describes the precise relationship.

\subsection{\boldmath Kernel of $\Psi:\Mod(S_{g,1})\to\Mod(X_{g}^{k})$}\label{sec:kernel}

Observe that the kernel of $\Psi$ is contained in the point-pushing subgroup $\pi_1(S_g)<\Mod(S_{g,1})$. This is because $\Psi$ composed with the natural map $\Mod(X_{g}^{k})\to\Mod(S_g)$ is the natural map $\Mod(S_{g,1})\to\Mod(S_g)$, whose kernel is the point-pushing subgroup. Thus we want to understand the image of the point-pushing subgroup under the section $\sigma$ used to define $\Psi$. What we find is a simple relationship between three surface group representations: 
\[\begin{xy}
(0,20)*+{\pi_1(S_g)}="A";
(-60,0)*+{\pi_1(S_g)}="B";
(60,0)*+{\pi_1(S_g)}="C";
(0,0)*+{\AUT\big(\pi_1(X_{g}^{k})\big)}="D";
{\ar "A";"D"}?*!/_25mm/{\text{\emph{inner auts of $\pi_1(S_g)$, lifted}}};
{\ar "B";"D"}?*!/^4mm/{\text{\emph{inner auts of $\pi_1(X_g^k)$}}};
{\ar "A";"D"}?*!/^3mm/{\sigma};
{\ar "C";"D"}?*!/_4mm/{\text{\emph{transvections}}}
\end{xy}\]

The main results are Proposition \ref{prop:section-on-inner} and Corollary \ref{cor:section-on-inner} below. In order to state Proposition \ref{prop:section-on-inner}, we need the following notation. Let \[\delta:H_1(S_g;\Z)\to H^1(S_g;\Z)\] be the Poincar\'e duality map, given explicitly by $\gamma\mapsto \langle -,\gamma\rangle$, where 
\[\langle -, -\rangle:H_1(S_g;\Z)\times H_1(S_g;\Z)\to\Z\] is the algebraic intersection form. We use $\hat\delta$ denote the composition
\[\hat\delta:H_1(S_g;\Z)\xrightarrow{\delta} H^1(S_g;\Z)\hookrightarrow\AUT\big(\pi_1(X_{g}^{k})\big).\]
This map is given explicitly by $\hat\delta(\gamma)(w)=w\cdot z^{\langle[\bar w],\gamma\rangle}$, where $\bar w$ is the image of $w$ under $\pi_1(X_{g}^{k})\to\pi_1(S_g)$ and $[\bar w]\in H_1(S_g;\Z)$ is the corresponding homology class. 

Fix a basepoint $\star\in S_{g,1}$. Recall that we have fixed a standard generating set $\{\alpha_i,\beta_i\}$ of $\pi_1(S_{g,1},\star)\cong F_{2g}$ so that $c:=\prod_i[\alpha_i,\beta_i]$ is a loop around the puncture $*$ of $S_{g,1}=S_g\setminus\{*\}$. 
Define 
\begin{equation}\label{eqn:change-basepoint}\Pi:\pi_1(S_{g,1},\star)\to\pi_1(S_g,*)\end{equation}
by $\gamma\mapsto\epsilon.\gamma.\bar\epsilon$, where $\epsilon$ is a fixed arc from $*$ to $\star$. 


\begin{proposition}\label{prop:section-on-inner}
Fix $t\in\pi_1(S_g,*)$, and let $\Push(t)\in\Mod(S_{g,1})\cong\Out^*(F_{2g})$ be the point-pushing mapping class. If $\tilde t\in\pi_1(S_{g,1},\star)$ is any lift of $t$ (i.e.\ $\Pi(\tilde t)=t$), then 
\begin{equation}\label{eqn:section-on-inner}\sigma\big(\Push(t)\big)=C_{\iota\tilde t}\circ \hat\delta([kt]),\end{equation}
Here $C_{x}$ denotes conjugation by $x$, and the maps $\iota:F_{2g}\to\pi_1(X_g^k)$ and $\sigma:\Out^*(F_{2g})\to \AUT\big(\pi_1(X_{g}^{k})\big)$ are defined in \eqref{eqn:iota} and \eqref{eqn:sigma}. 
\end{proposition}

As a sanity check, observe that $C_{\iota\tilde t}$ does not depend on the choice of lift $\tilde t$ because any two lifts differ by an element of the normal closure of $c$ in $\pi_1(S_{g,1},\star)=F_{2g}$, and conjugation by any such element is trivial on $\pi_1(X_{g}^{k})$. 

\begin{proof}[Proof of Proposition \ref{prop:section-on-inner}]

It suffices to prove the lemma for $t\in\pi_1(S_g,*)$ that are represented by a non-separating simple closed curve. To see this, first note that $\pi_1(S_g,*)$ is generated by these curves. Furthermore, the groups $\Inn\big(\pi_1(X_{g}^{k})\big)$ and $H^1(S_g;\Z)$ commute in $\Aut\big(\pi_1(X_{g}^{k})\big)$, so
\[\big[C_{\iota \tilde t_1}\circ\hat\delta([t_1])\big]\circ\big[C_{\iota\tilde t_2}\circ\hat\delta([t_2])\big]=C_{\iota(\tilde t_1*\tilde t_2)}\circ \hat\delta([t_1*t_2]).\]

Assume now that $t\in\pi_1(S_g,*)$ is represented by a non-separating simple closed curve. After an isotopy, we can assume that $t$ contains $\epsilon$ as a sub-arc. Choose $\tilde t$ as pictured in Figure \ref{fig:rep}.
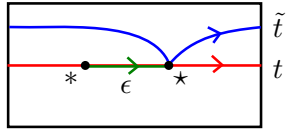
\begin{figure}[h!]
\labellist
\pinlabel $*$ at 230 430
\pinlabel $\star$ at 290 430
\pinlabel $t$ at 343 435
\pinlabel $\tilde t$ at 343 459
\pinlabel $\epsilon$ at 260 425
\endlabellist
\centering
\includegraphics[scale=.7]{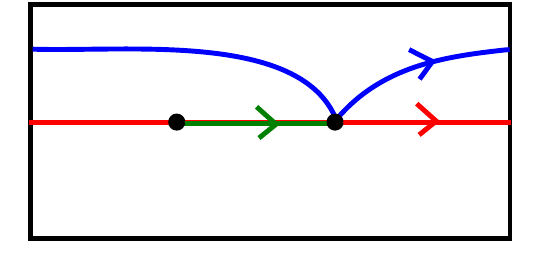}
\caption{A small regular neighborhood of a loop representing $t\in\pi_1(S_g,*)$ and a lift $\tilde t\in\pi_1(S_{g,1},\star)$.}
\label{fig:rep}
\end{figure}

We want to show that 
\[\sigma\big(\Push(t)\big)(w)=\big[C_{\iota\tilde t}\circ\hat\delta([t])\big](w)\]
for each $w\in\pi_1(X_{g}^{k})$. Since this is obviously true for $w=z$, it suffices to show this equality for $w=\iota(s)$ for $s\in \pi_1(S_{g,1},\star)$; furthermore, it suffices to show the equality on any generating set of $\pi_1(S_{g,1},\star)$. We use the (infinite) generating set consisting of curves of one of the forms pictured in Figure \ref{fig:gens} (the intersection of these curves with the annulus around $t$ has one component). 

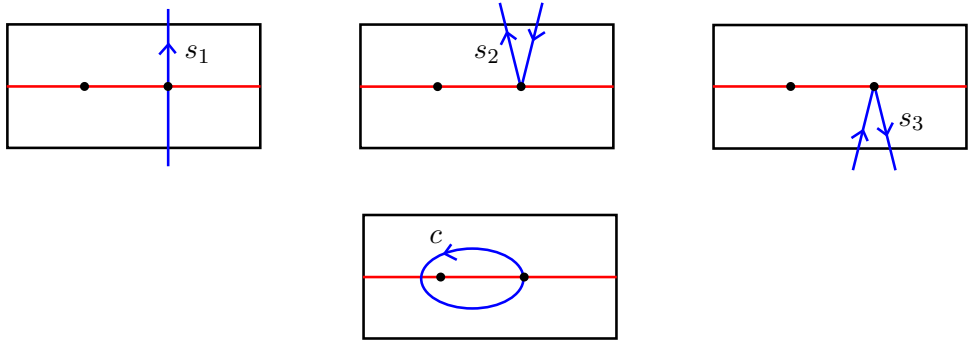
\begin{figure}[h!]
\labellist
\pinlabel $c$ at 235 237
\pinlabel $s_1$ at 105 337
\pinlabel $s_2$ at 263 337
\pinlabel $s_3$ at 495 300
\endlabellist
\centering
\includegraphics[scale=.7]{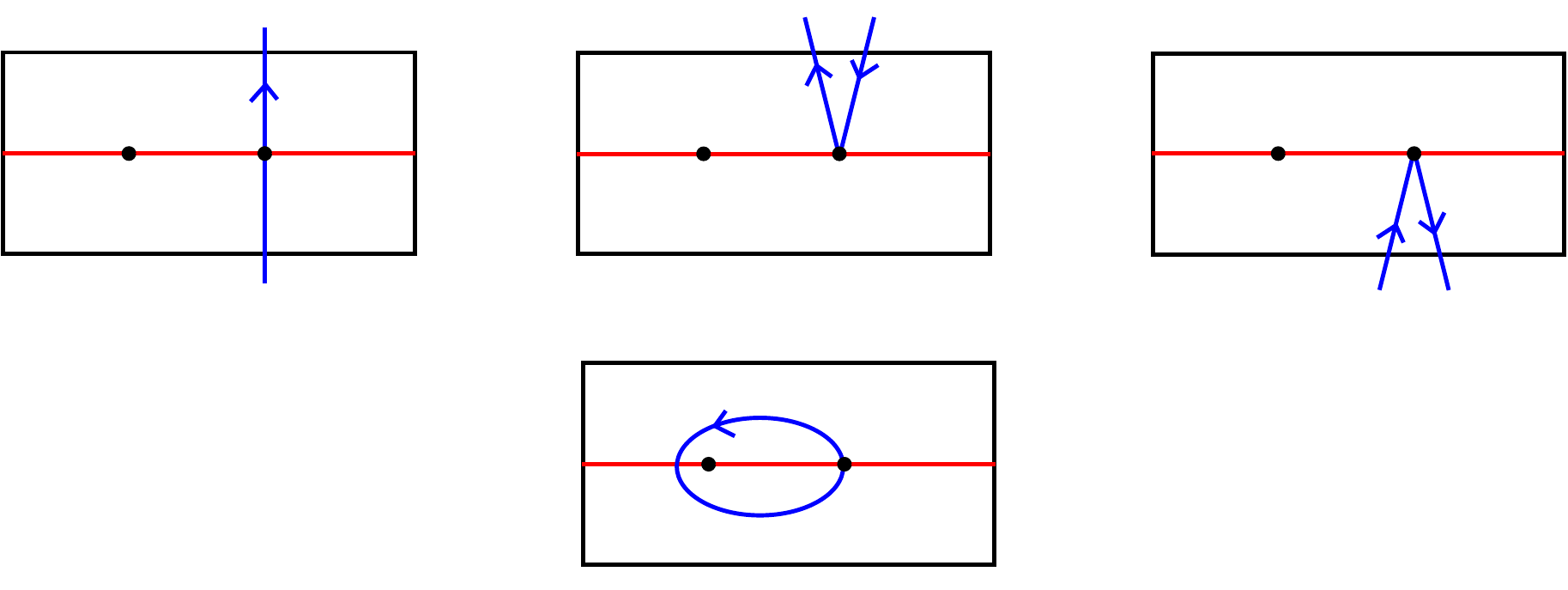}
\caption{The group $\pi_1(S_{g,1},\star)$ is generated by $\tilde t$ and loops of the form pictured above. }
\label{fig:gens}
\end{figure}

Note that $\Push(t)$ fixes the basepoint $\star$, so we can compute the action of $\Push(t)$ on $s\in\pi_1(S_{g,1},\star)$. We compute the action of $\Push(t)$ on the elements in Figure \ref{fig:gens} as follows. See Figure \ref{fig:action} for an illustration. 
\[s_1\mapsto (\tilde t)^{-1}s_1\tilde t c^{-1}\>\>\>\text{ and }\>\>\>s_2\mapsto (\tilde t)^{-1}s_2\tilde t\>\>\>\text{ and }\>\>\>s_3\mapsto c(\tilde t)^{-1}s_3\tilde tc^{-1}\>\>\>\text{ and }\>\>\>c\mapsto c.\]
\vspace{.1in}

\begin{figure}[h!]
\labellist
\pinlabel $\Push(t)(s_1)$ at 100 100
\pinlabel $\Push(t)(s_2)$ at 263 100
\pinlabel $\Push(t)(s_3)$ at 465 60
\endlabellist
\centering
\includegraphics[scale=.7]{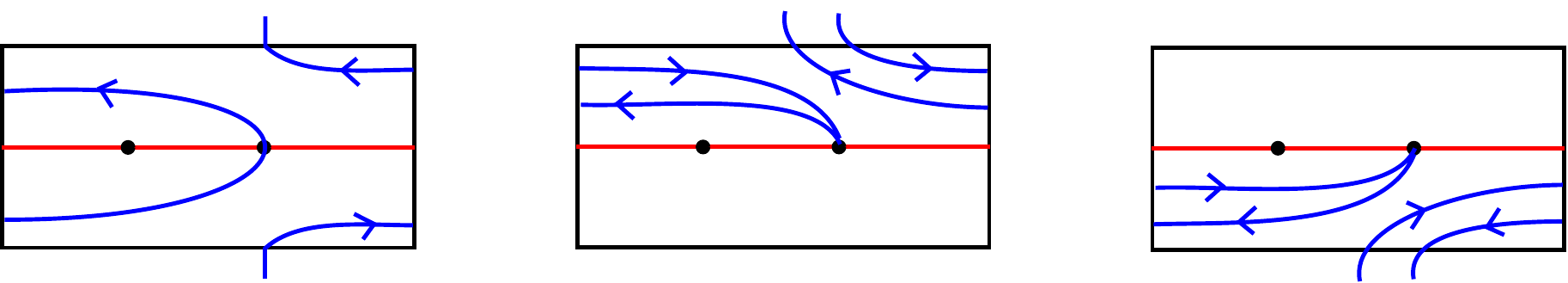}
\caption{Action of point-pushing about $t$ on the loops in Figure \ref{fig:gens}. The curve $c$ is fixed up to isotopy (up to isotopy $\Push(t)$ is the identity on a neighborhood of $t$ that contains $c$). }
\label{fig:action}
\end{figure}

This proves that, for example, that 
\[\sigma\big(\Push(t)\big)(\iota s_1)= (\iota \tilde t)^{-1}(\iota s_1)(\iota \tilde t) z^{-k}=\big[C_{\iota\tilde t}\circ \hat\delta([kt])\big](\iota s_1).\]
We conclude similarly for the generators $s_2,s_3$. This proves the desired formula for $\sigma\big(\Push(t)\big)$. 
\end{proof}

The following corollary is an immediate consequence of Proposition \ref{prop:section-on-inner}. 
\begin{corollary}\label{cor:section-on-inner}
Consider the composition 
\begin{equation}\label{eqn:section-to-Out}\Psi:\AUT\big(\pi_1(S_g)\big)\xrightarrow{\sigma}\AUT\big(\pi_1(X_{g}^{k})\big)\to\OUT\big(\pi_1(X_{g}^{k})\big).\end{equation}
The restriction of $\Psi$ to $\pi_1(S_g)\cong\Inn\big(\pi_1(S_g)\big)$ factors as follows. 
\[\begin{xy}
(-20,0)*+{\pi_1(S_g)}="A";
(40,0)*+{\AUT\big(\pi_1(S_g)\big)}="B";
(-20,-15)*+{H^1(S_g;\Z)}="C";
(40,-15)*+{\OUT\big(\pi_1(X_{g}^{k})\big)}="D";
{\ar "A";"B"}?*!/_3mm/{\text{\emph{conjugation}}};
{\ar "B";"D"}?*!/_3mm/{\Psi};
{\ar "A";"C"}?*!/^6mm/{k\delta\circ\ab};
{\ar "C";"D"}?*!/^3mm/{\text{\emph{transvections}}}
\end{xy}\]
Here $\ab$ denotes the abelianization map $\pi_1(S_g,*)\to H_1(S_g;\Z)$. 
\end{corollary}

\subsection{Proof of Theorem \ref{thm:main}} 

Using the isomorphisms between mapping class groups and automorphism groups, the desired diagram is equivalent to the following one.
\[\begin{xy}
(-60,0)*+{1}="A";
(-40,0)*+{\pi_1(S_g)^{\ab}}="B";
(0,0)*+{\AUT\big(\pi_1(S_g)\big)/\pi'}="C";
(40,0)*+{\OUT\big(\pi_1(S_g)\big)}="D";
(60,0)*+{1}="E";
(-60,-20)*+{1}="F";
(-40,-20)*+{\Hom(\pi_1(S_g),\Z)}="G";
(0,-20)*+{\OUT\big(\pi_1(X_{g}^{k})\big)}="H";
(40,-20)*+{\OUT\big(\pi_1(S_g\big))}="I";
(60,-20)*+{1}="J";
{\ar "A";"B"}?*!/_3mm/{};
{\ar "B";"C"}?*!/_3mm/{};
{\ar "C";"D"}?*!/^3mm/{};
{\ar "D";"E"}?*!/^3mm/{};
{\ar "F";"G"}?*!/_3mm/{};
{\ar "G";"H"}?*!/_3mm/{};
{\ar "H";"I"}?*!/^3mm/{};
{\ar "I";"J"}?*!/^3mm/{};
{\ar "B";"G"}?*!/_3mm/{k\delta};
{\ar "C";"H"}?*!/^3mm/{};
{\ar@{=} "D";"I"}?*!/_3mm/{}
\end{xy}\]

The map $\Psi$ in \eqref{eqn:section-to-Out} descends to the middle vertical map and restricts to the left vertical map by Corollary \ref{cor:section-on-inner}. The fact that $\sigma$ is a section (Corollary \ref{cor:Aut-SES}) implies that the middle vertical map descends to the identity map on $\OUT\big(\pi_1(S_g)\big)$. When $k=1$, the middle vertical map is an isomorphism by the five lemma. This concludes the proof of Theorem \ref{thm:main}.

\section{Spectral sequence computation}\label{sec:sseq}

In this section we prove Theorem \ref{thm:H2H1}. This is achieved by two different computations using the Lyndon--Hochschild--Serre (LHS) spectral sequence. Recall that this spectral sequence takes input a short exact sequence of groups $1\to N\to G\to Q\to 1$ and a $G$-module $A$, has $E_2$ page 
\[E_2^{p,q}=H^p\big(Q;H^q(N;A)\big),\]
and converges to $H^{p+q}(G;A)$. For both computations we use the Birman exact sequence, but with different choices of the module $A$. 

{\it Notational note.} To simplify the notation, we use the convention that cohomology groups have $\Z$ coefficients unless otherwise specified.

\subsection{Euler class computation}

Our goal in this section is to prove Proposition \ref{prop:euler} below, which implies Corollary \ref{cor:splitting}. 

\begin{proposition}\label{prop:euler}
Fix $g\ge1$. Let $eu_k$ be the Euler class of the extension \eqref{eqn:Mod-SES}. Then $eu_k=k\> eu_1$, and $eu_1$ has order $2g-2$ in $H^2\big(\Mod(S_g);H^1(S_g)\big)$. 
\end{proposition}

\begin{proof}
The relation $eu_k=k\>eu_1$ already follows from Theorem \ref{thm:main}. Indeed, choosing a set-theoretic section for the sequence in the top row of the diagram in Theorem \ref{thm:main} gives a cocycle representative for $eu_k$ that is $k$ times the cocycle representative for $e_1$. 

Now we prove that $eu_1$ generates a cyclic subgroup isomorphic to $\Z/(2g-2)\Z$ in $H^2\big(\Mod(S_g);H^1(S_g)\big)$. Our method is to apply the LHS spectral sequence to the Birman exact sequence with the module $A=H^1(S_g)$. Here
\[E_2^{p,q}\cong H^p\big(\Mod(S_g);H^q(S_g;A)\big).\]
A portion of the associated 5-term exact sequence is as follows.
\[\begin{array}{rll}0\to H^1\big(\Mod(S_g);H^1(S_g)\big)\to H^1\big(\Mod(S_{g,1});H^1(S_g)\big)&\xrightarrow{A} \Hom\big(H_1(S_g), H^1(S_g)\big)^{\Mod(S_g)}\\[2mm]
&\xrightarrow{d_2^{0,1}}H^2\big(\Mod(S_g);H^1(S_g)\big)
\end{array}
\]
This sequence has been studied by Morita. Morita \cite[Prop.\ 4.1]{Morita} computes that the first term vanishes, so the map $A$ is injective. The group $\Hom\big(H_1(S_g), H^1(S_g)\big)^{\Mod(S_g)}$ is isomorphic to $\Z$ and generated the Poincar\'e duality isomorphism $\delta$. Morita \cite[proof of Prop.\ 6.4]{Morita} shows that the image of $A$ is $(2g-2)\Z$. Consequently, the differential $d_2^{0,1}$ descends to an injection $\Z/(2g-2)\Z\to H^2\big(\Mod(S_g);H^1(S_g)\big)$. 

It remains to show that $d_2^{0,1}$ sends a generator to $eu_1$. 
The differential $d_2^{0,1}$ is the transgression; see e.g.\ \cite[Prop.\ 1.6.6, Thm.\ 2.4.3]{Neukirch}. By standard knowledge of the transgression applied to our situation, we find that $d_2^{0,1}$ sends a generator to $\delta_*(eu)$, where $eu$ is the Euler class of the extension (\ref{eqn:BES-quotient}), and 
\[\delta_*: H^2\big(\Mod(S_g);H_1(S_g)\big)\to H^2\big(\Mod(S_g);H^1(S_g)\big)\]
is the isomorphism induced by the Poincar\'e duality isomorphism $\delta$. (For this property of the transgression, see \cite[\S I.6, Exercise 1-2]{Neukirch}. While that reference is mainly concerned with finite or profinite groups, the analysis of the transgression contained given there applies more generally.) Finally, we observe that $\delta_*(eu)=eu_1$ by Theorem \ref{thm:main}. 
\end{proof}

\subsection{Computation of $H^2\big(\Mod(S_g);H^1(S_g)\big)$}

Running the LHS spectral sequence with the trivial module $A=\Z$, we prove that if $g\ge8$, then 
\begin{equation}\label{eqn:H2H1}H^2\big(\Mod(S_g);H^1(S_g)\big)\cong\Z/(2g-2)\Z.\end{equation}
Combining this with Proposition \ref{prop:euler} proves Theorem \ref{thm:H2H1}. 
The relevant portion of the spectral sequence appears below. 


\begin{tikzpicture}
  \matrix (m) [matrix of math nodes,
    nodes in empty cells,nodes={minimum width=3ex,
    minimum height=5ex,outer sep=0pt},
    column sep=0.1ex,row sep=1ex]{
                  &                                                &      &     & \\
          2      &   H^0(\Mod(S_g);H^2(S_g))   &     &     & \\
          1      &  0                                           &  0  & H^2(\Mod(S_g);H^1(S_g)) & \\
          0      &  \Z                                          & 0 &H^2(\Mod(S_g))  & H^3(\Mod(S_g)) & H^4(\Mod(S_g)) \\
                  &   0                                          &  1  &    2  & 3 &4 \\};
\draw[->] (m-2-2) -- (m-3-4) node[midway, above] {$d_2^{0,2}$};
  \draw [->] (m-3-4) -- (m-4-6) node[midway, above] {$d_2^{2,1}$};
\draw[thick] (-6.5,-2.1) -- (-6.5,2) ;
\draw[thick] (-7.2,-1.7) -- (7.3,-1.7) ;
\end{tikzpicture}

%

The computations in the first column are easy. In the second column, Morita \cite[Prop.\ 4.1]{Morita} computed $H^1\big(\Mod(S_g);H^1(S_g)\big)=0$ for $g\ge1$. The other computation $H^1\big(\Mod(S_g)\big)=0$ holds for $g\ge1$ because the abelianization of $\Mod(S_g)$ is finite \cite[\S5.1.2-3]{Farb-Margalit}.  

According to \cite[Cor.\ 1.2]{Tillmann}, 
\[H_*\big(\Mod(S_{g,1})\big)\cong H_*\big(\Mod(S_g)\big)\otimes \Z[x]\]
in degrees $g\ge 2*$. Here $x$ has degree 2. 
Applying this and using the universal coefficients theorem, we conclude that 
\[H^i\big(\Mod(S_g)\big)\to H^i\big(\Mod(S_{g,1})\big)\]
is an isomorphism if $i=3$ and $g\ge6$, and it is injective if $i=4$ and if $g\ge8$. 

Since the map $H^4\big(\Mod(S_g)\big)\to H^4\big(\Mod(S_{g,1})\big)$ is injective, the differential $d_2^{2,1}$ is zero. Since the map
$H^3\big(\Mod(S_g)\big)\to H^3\big(\Mod(S_{g,1})\big)$ is an isomorphism, the differential $d_2^{0,2}$ is surjective. 

Thus, the filtration of $H^2\big(\Mod(S_{g,1})\big)$ coming from the $E_\infty$ page gives an exact sequence 
\[\begin{array}{rcl}0\to H^2\big(\Mod(S_g)\big)\to H^2\big(\Mod(S_{g,1})\big)&\xrightarrow{F}&H^0\big(\Mod(S_g);H^2(S_g)\big)\cong\Z\\
&\xrightarrow{d_2^{0,2}}& H^2\big(\Mod(S_g);H^1(S_g)\big)\to0.\end{array}\]
For $g\ge4$, 
\[H^2\big(\Mod(S_g)\big)\cong\Z[e_1]\>\>\>\text{ and }\>\>\>H^2\big(\Mod(S_{g,1})\big)\cong\Z[e,e_1]\]
and the map $\Z[e_1]\to\Z[e,e_1]$ is the obvious one $e_1\mapsto e_1$. We claim that $F(e)=2-2g$. From this we deduce the desired isomorphism (\ref{eqn:H2H1}). The claim follows from the fact that the extension that defines $e$, when restricted to the point-pushing subgroup $\pi_1(S_g)<\Mod(S_{g,1})$, gives the extension 
\[1\to\Z\to\pi_1(US_g)\to\pi_1(S_g)\to1\]
where $US_g$ is the unit tangent bundle. See \cite[\S5.5.5]{Farb-Margalit}. This extension has Euler class $2-2g$, so the claim follows.

\bibliographystyle{alpha}
\bibliography{citing}

\end{document}